\documentclass[a4paper,10pt]{article}
\usepackage{amssymb}
\usepackage[T2A]{fontenc}
\usepackage[cp1251]{inputenc}
\usepackage{amsmath}

\usepackage[dvips]{graphicx}
\usepackage[english,russian]{babel}

\newtheorem{theorem}{Теорема}[section]

\newtheorem{definitionhead}[theorem]{Определение}

\author{А.~К.~Ковальджи, А.~Я.~Канель-Белов}
\title{Занятия по математике -- листки и диалог}

\begin{document}
\maketitle




{\small
\begin{quote}
Однажды я провел в размышлениях целый день без еды и целую ночь без сна, но я ничего не добился. Было бы лучше посвятить то время учению.
Тот, кто учится не размышляя, впадет в заблуждение. Тот, кто размышляет, не желая учиться, окажется в затруднении.

Учиться и не думать -- бесполезно, а думать и не учиться -- опасно.

\begin{flushright}
Конфуций
\end{flushright}

\end{quote}

}

\section{Введение}
 В кружках и математических классах получила распространение т.н. ``листковая система’’, которая в силу технического удобства повсеместно вытесняет иные формы занятий.\footnote{Похожая ситуация сложилась с массовым тестированием (IQ, ЕГЭ, ``Кенгуру").}

Формы листковой системы могут варьироваться, но, в общем и целом, она заключается в следующем:

\begin{enumerate}
 \item Учащиеся получают
 листки с задачами. Иногда туда включается минимум теории.
 \item Решив задачу, учащийся поднимает руку, к нему подходит проверяющий.
 \item Если задача решена правильно, учащемуся в ведомость вносится знак ``плюс''.
 \item Иначе ему приходится думать дальше.
 \item Имеются требования по количеству решенных задач к определенному сроку с учетом их сложности, некоторые задачи объявляются обязательными для решения.
\end{enumerate}

Одна из целей листков –– освоение новой области математики в процессе работы. Сначала проводится ``ликбез'', когда перечисляются необходимые факты и понятия, а дальше, по мере возрастания трудности задач, ученики самостоятельно выстраивают необходимую теорию, и роль учителя заключается в проверке и корректировке их самостоятельной работы.

 Система листков сыграла позитивную роль и в руках профессионалов дала прекрасные результаты, но со временем она стала повсеместно использоваться и во многом потеряла изначальные смыслы и традиции. Поэтому возникла необходимость обсудить педагогические цели, которым она служит, а также иные цели, для достижения которых требуется диалог и другие формы работы, незаслуженно отодвинутые в сторону.

Мы обсуждаем историю вопроса, затем -- достоинства и недостатки листковой системы. Далее приводим примеры тем занятий, где необходим диалог. В конце мы говорим о синтезе системы листков и диалога.

\section{История вопроса}

Современный тип математического кружка для школьников является гениальным изобретением Д.~О.~Шклярского и есть следствие сделанных им нововведений.
\footnote{К сожалению, до сих пор не собраны воспоминания о работе Д.~О.~Шкляр\-ско\-го и его педагогической технологии.} По всей видимости, технологию Шклярского вызвали к жизни математические олимпиады, получившие тогда распространение.
\footnote{Работы Г.~Полиа и Г.~Сеге, а также венгерская комбинаторика (знаменитая т.н. ``венгерская математика'') вызваны к жизни олимпиадами.} Кружок, основанный на решении задач, стал достаточно эффективным и дал замечательные результаты.

Вместе с тем обучение некоторым сторонам деятельности математика неудобно проводить в формате олимпиадного кружка, так что актуальным стал поиск дополнительных форм. Например, проектный (учебно-исследовательский) подход. Научные конференции школьников зачастую проводятся сразу по нескольким предметам, что создает условия для общения школьникам с разной специализацией. Самостоятельные исследования являются альтернативной формой деятельности, позволяющей увидеть и развить у школьника  полезные для научной работы качества, которые олимпиады не раскрывают. Эти подходы взаимно дополнительны. Надо стремиться к тому, чтобы преподаватели, исповедующие тот или иной подход, не конкурировали между собой, а сотрудничали, понимали друг друга. Проектному подходу посвящены работы \cite{Conferencii,Roitberg,Sgibnev}.

В середине 60-х -- конце 70-х олимпиадное движение достигло определенного пика развития. В жюри олимпиад активно работали такие математики, как А.~Н.~Колмогоров, И.~М.~Гельфанд, Е.~Б.~Дынкин, В.~И.~Арнольд, Н.~Б.~Васильев и другие, многое в технологии проведения олимпиад было наработано (хотя далеко не всё. Не осознавалась роль утешительных задач, варианты были хуже сбалансированы, чем современные, однако вкус и научное содержание были в целом выше). Развившийся задачный подход вызвал к жизни попытки разбиения теоретического материала на задачи.
Элементы ``листковой системы'' возникли в 60-е годы, в частности, в Вечерней математической школе при Московском математическом обществе,
 созданной Е.~Б.~Дынкиным.
\footnote{Авторам представляется крайне важным собрать воспоминания о методике работы этого выдающегося математика, педагога и организатора образования.}

\subsection{Появление листков}
Чуть позднее усилиями Н.~Н.~Константинова в Москве сложилась форма обучения в кружках по математике -- листки с задачами, которые выдавались ученикам на занятиях. \footnote{Затем в Ленинграде независимо появилась похожая система, разработанная В.~П.~Федотовым. Ленинградская (петербургская) система преподавания, созданная в настоящем виде С.~Е.~Рукшиным, нуждается в отдельной статье.} Каждый ученик решал эти задачи в индивидуальном темпе, а учитель проверял правильность решений, делал замечания и давал советы. Листки, как правило, бывают тематическими и рассчитаны на определенный возраст и уровень подготовки учеников. Эти листки стали накопителями идей и популярных задач, которые передавались по городам и весям, помогая всем учителям ценными наработками талантливых математиков-педагогов.

Система листков в 57-й московской школе подробно описана в \cite{DavidovichPushkarChekanov,MatemvZadachah}.
При ее осуществлении сложились традиции, связанные в частности, с большим числом проверяющих.

Она позволяет начинающему преподавателю начать работать. Студенты не только помогают проводить занятие, но и образуют промежуточное звено между старшим преподавателем и учеником. Так, о молодых преподавателях в книге \cite{MatemvZadachah} написано: ``Они лучше чувствуют ребенка -- между ними нет психологического барьера (и потому неудивительно, что общение школьников и студентов не ограничивается рамками школьных уроков -- это и походы, песни под гитару, обсуждение книг и фильмов; причем все это часто продолжается и после выпуска). У них огромное желание поделиться тем, чему их самих научили в школе и в вузе. Наконец, они еще помнят, как их учили; причем не только то, чт\'о получалось, но и то, чт\'о преподаватели по их (выпускников) мнению делали неудачно. Поэтому им практически и не нужно специальное педагогическое образование -- они сразу готовы учить по данной системе, разумеется, при чутком руководстве''. Кроме того, листковая система в 57-й школе не так примитивна, как это зачастую наблюдается в иных местах, поскольку люди осознают и решают непростые педагогические задачи (см. предисловие к книге \cite{MatemvZadachah}).

В то же время, как писал И.~С.~Рубанов, ``работа по листочкам очень заманчива: преподавателям не надо диктовать, детям –- записывать, материалы занятий остаются в упорядоченном виде. Но она сопряжена с несколькими опасностями. Первая: тематический листочек сам по себе –- сильная подсказка. Вторая –– листочек сковывает преподавателя, лишает гибкости, сильно ограничивает возможность импровизации. Третья: листочек, раскрывая наперед все карты, лишает занятие интриги. Выход может быть в частичном отказе от листочков, выдаче их не в начале занятия, дроблении на фрагменты, выдаваемые в нужные моменты. Может быть, иногда вообще имеет смысл распечатывать задачи и теоретические комментарии поодиночке, устраивая каждому ученику индивидуальную траекторию?''

При дальнейшем распространении листковой системы произошло упрощение процесса. При этом достоинства (описанные в \cite{DavidovichPushkarChekanov,MatemvZadachah}) уменьшились, а недостатки усилились, особенно в ситуации одного-двух педагогов на класс и недостаточных традиций. Ослабление методических требований к преподавателям
в сочетании с появлением тренеров, специализирующихся на подготовке к олимпиадам (и далеких от научного исследования) сформировало специфический слой людей. В олимпиадном движении накопились проблемы, см. \cite{dver}. Эти проблемы усугубили недостатки листковой системы, особенно в руках тренеров. Упростилось и ухудшилось преподавание в некоторых лагерях. Технологизм оказался соблазнительным.

Система рейтингов не обязательна при листковой системе. Но она возникает естественным образом, хотя и не всегда. Вот отзыв участника: ``Спецификой листковой и рейтинговой системы является дополнительный стимул оказаться быстрее, выше и сильнее, сдать много задач. Это мотивирует постоянно решать задачи. В то же время положение в рейтинге может демотивировать слабую половину группы. Элемент соревнования может негативно влиять на коллективизм.''
\footnote{К сожалению, рейтингами увлекаются отнюдь не только подростки, но и власть предержащие, попадая при этом в ловушку подчинения сомнительным рейтинговым агентствам. Упомянем нелепую программу ``топ-100'' -- добиться включения нескольких российских вузов в 100 вузов наибольшего рейтинга. Согласно результатам престижной международной студенческой олимпиады http://www.imc-math.org/ , мехмат МГУ и Физтех обычно входят в пятерку лучших. При этом Физтех до недавнего времени имел 400-й рейтинг. Кроме того, около $20\%$ выпускников мехмата уже имеют научные результаты, пригодные к публикации.
Показатель, неслыханный в западном мире. В то же время рейтинг мехмата довольно низкий. Такой подход к оценке образования антинаучен -- следует изучать взаимодействие студента с вузом, принимая во внимание особенности, как студента, так и страны, в которой он обучается. Только так можно заимствовать то лучшее, что есть на ``Западе'' (кстати, ``Запад'' -- это отнюдь не одна точка, а страны с разными культурными традициями). О принципиальной порочности рейтинговой системы см. http://www.mccme.ru/free-books/bibliometric.pdf.}

Примечателен отзыв другого бывшего олимпиадника: ``Начальство судит о крутости кружка по успехам на соревнованиях. Школьник идёт на кружок, если там учат побеждать на соревнованиях. Так было и будет. '' Листковая система иногда преподносится, как единственно возможная. Вот реакция современного олимпиадника на обсуждение иных форм работы: `` Сложилось ощущение, что статья написана для кружка, где есть сильный препод или сильный ученик. Если нет ни того, ни другого, то отличные от листковой системы методы надо использовать только в виде исключения''.

\section{Система листков}

Фактически -- это система самообучения ученика, предполагающая его инициативу под руководством учителя. Для сильных учеников, умеющих самостоятельно работать, занятие по листкам эффективно, а начинающих кружковцев листки могут разочаровать.\footnote{Как сообщает С.~В.~Петухов: Так получилось, что, в основном, я работаю, как раз, со ``средними
детьми''. Им, действительно, сложно и не всегда интересно работать по листочку. Действительно, добавляю элементы диалога, разбиваю листочек на 3--4
 части -- и дело идёт намного лучше. Даже более того, листочек без изменений
 и дополнений считаю неэффективным, т.к. крайне мало детей способны его
 разумно решать без дополнительной мотивации и объяснений.}

Данная система эффективно решает следующие педагогические задачи.

\begin{enumerate}
 \item Начинающий учитель (даже старшеклассник), получив набор листков, может сразу начать проводить занятия. При этом, конечно, его квалификации должно хватить на то, чтобы проверить правильность логики ученика, который придумал другое решение задачи. (В дальнейшем, однако, рост преподавателя затрудняется, поскольку не возникает остро ощущаемой потребности в умении выступать перед классом и самому готовить занятия).
 \item Общение учеников с проверяющими близкого возраста.

 \item Индивидуальный подход к учащемуся. Возможность решать задачи в индивидуальном темпе.

 \item Удобство технической ``дрессировки'' и проработки малоинтересных, но необходимых деталей.

 \item Реализация  подхода, при котором освоение материала во многом сводится к последовательности решенных задач.
 \item Возможность лучше освоить тему самому без посторонней помощи.
 Доказательство, придуманное  самим учеником, дорогого стоит.
 \item Листки с предисловием позволяют разобраться в вопросе на хорошем уровне, что используется в дистанционном образовании, например, в ВЗМШ
(Всероссийской заочной многопредметной школе),)
созданной И.~М.~Гельфандом.

\end{enumerate}

В то же время продолжением достоинств являются недостатки. Критический разбор отнюдь не означает ``борьбы'' с листковой системой. Наша цель -- осознать ее ограничения и призвать к необходимости ее дополнения иными формами работы.

\begin{enumerate}
 \item Преподаватель довольно быстро входит в курс дела, но при этом не стимулирован к росту. В результате недостаточно понимается важность педагогической подготовки преподавателя (а лучше сказать –- его психологической готовности к работе с детьми).
 \item Ученики не видят разные решения, не учатся выступать, слушать друг друга, спорить. Между тем все это так же необходимо, как и умение решать задачи.

 \item Не стимулируется самостоятельный творческий поиск преподавателя, который приучается работать по готовым чужим материалам.
( На Малом мехмате
в 90-е годы был возмутительный случай, когда старший по параллели выгнал студента, который вел занятия не по централизованно раздаваемым листкам.)
 \item Реализуется узко технологический подход к математике, во многом эффективный. Однако он всё подминает под себя, создавая узость математического видения. Математика обсуждается локально, тактически в рамках частной задачи, а к стратегическому \footnote{Т.е.
на уровне плана решения, осознания причин успеха, движущих идей.} обсуждению листковый метод приспособлен плохо. За деревьями ученик должен видеть лес.
 \item В рамках исключительно листкового подхода затрудняется и ограничивается развитие теоретического мышления, о чем мы поговорим ниже.

 \item Преподавателю скучно повторять одно и то же разным ученикам в течение занятия, ему приходится по многу раз объяснять одни и те же ошибки.
 \item Слабый учащийся или начинающий может ничего не решить, в то же время проверяющие зачастую уделяют больше времени сильным учащимся. Отметим, что учащиеся {\it растут неравномерно}. Слабый сегодня может оказаться сильным завтра, если только не проявлять к нему безразличие, что зачастую допускают молодые преподаватели. Ситуация несколько улучшается при наличии утешительных задач.

 \item Нужно не только учиться все делать самому, но и уметь учиться у других, уметь работать
в команде. Однако листковая система не способствует этому.
\end{enumerate}

\subsection{Связь системы листков с олимпиадами}

Ситуация усугубляется тем, что математическое образование зачастую сводят к подготовке к олимпиадам. Так получилось, что б\'ольшая часть внеклассной работы школьников прямо или косвенно связана с олимпиадами или иными соревнованиями по решению четко поставленных задач в жестких временных рамках. Даже если руководитель кружка заявляет, что не занимается подготовкой к олимпиадам, материалы кружка в сильной степени основаны на олимпиадных задачах. Олимпиадный подход имеет свои достоинства и недостатки, ему посвящена статья \cite{dver}. При этом преподаватель зачастую попадает в ловушку технологизма и утилитаризма, когда достижение локальных целей заслоняет более высокие ценности. Наша критика относится не столько к самой листковой системе, сколько к ее роли в сложившейся ситуации.

\subsection{Теоретическое и стратегическое мышление}

Математика --- это не только наука о решении
конкретных
задач. Это и взгляд на мир, и философия. Математик
также
строит новые теории, вводит и осмысляет новые понятия. Нужно уметь не только решать задачи, но и их ставить, думать о {\it мотивировках}.
Наряду с {\it тактикой} есть
{\it стратегия} и даже {\it надстратегия} (иногда ее называют ``философией'' или ``идеологией'').
\footnote{Примечательно, что один олимпиадный деятель оспаривал само наличие стратегического и тактического мышления. Между тем
помимо преодоления технических трудностей бывает необходимо рассуждать на макроуровне, посмотреть на ситуацию в целом,
понять смысл задачи.} Помимо тренировки в решении конкретных задач, важно развивать {\it теоретическое мышление}. Кстати сказать, оно чрезвычайно полезно и для олимпиадных успехов. Бывают
задачи, которые не решить напрямую. Надо {\it понять смысл} условия задачи (см. обсуждение задачи про ломаную, делящую квадрат на две равные части, в \cite{dver}, см. также \cite{BelovMarkelov}), сменить подход и исследовать ситуацию, создав (мини или не мини) теорию. Кстати, работа над достаточно содержательной задачей отнюдь не заканчивается
вместе с её решением,
см. \cite{Akvarium}. (Правда, ``философствование'' рискует оказаться пустым, и отчасти задачный подход есть реакция на это.)

Например, доказательство на языке $\varepsilon - \delta$ того, что $\lim\limits_{x\to 1}\frac{x^2-1}{x-1}=2$, отнюдь не служит удачным образцом воспитания теоретического мышления.
\footnote{При выборе дополнительных курсов для школьников зачастую идут по пути копирования вуза, в частности в преподавании математического анализа. Как следствие получается, что бывшие матшкольники на первых курсах бездельничают, а
иногда -- разучаются работать. У них возникают проблемы, вплоть до исключения. Выбору курсов для школьников следовало бы посвятить отдельную статью, здесь же мы ограничиваемся указанием на проблему.
Матанализ
может быть сервисом для иных курсов, в особенности физики, и чрезвычайно полезным может быть подход Я.~Б.~Зельдовича
(см. \cite{Zeld1,Zeld2}). Его
 весьма эмоционально --- и не всегда корректно --- оспаривал Л.~С.~Понтрягин, но затем с ним согласился.

 Использование же аксиоматики вещественного числа как полигона для воспитания культуры мышления представляется авторам неправильным. Прежде всего, теория действительных чисел становится интересной только при обсуждении более продвинутых вопросов, а отнюдь не на начальной стадии, да и сама идеология Дедекинда --- Вейерштрасса устарела после развития теории моделей -- более правильно развивать мышление в курсе математической логики, доступной для школьников, умеющих программировать. Есть и другие весьма интересные сюжеты, не отраженные в вузовском преподавании. Не случайно создатели Второй школы -- выдающиеся ученые -- использовали термин ``спецматематика''.}

Не нужно доказывать эквивалентность трех определений комплексного числа, а полезно обсудить разные интерпретации комплексных чисел. Вместе с тем,
имеет смысл анализировать разные определения выпуклости фигур, поскольку это дается сравнительно легко и помогает решать задачи. А в ``игре''
с определениями комплексных чисел, с одной стороны, имеется стремление подражать ``большой науке'', а с другой~-- мало содержания.

\section{Диалог}

Работу со школьниками один из авторов начал в 1973 году во Всесоюзной заочной математической школе,
где было много хороших задач. Самое интересное для него –– возникла переписка с талантливыми учениками, например, с Александром Гончаровым из Никополя, ныне известным математиком. В переписке обсуждались обобщения задач, родственные задачи, различные подходы к одним и тем же задачам.

Такой диалог преподавателя с учеником очень важен. В работе кружка
первично содержание занятий, а не их форма, но чем младше школьники, тем форма для них важнее. Многие дети еще не разобрались, нравится им математика или нет, поэтому они больше стремятся к развлечению и общему развитию, чем к знаниям и математической культуре. Задача учителя –– раскрепостить детей, разрешить им задавать много вопросов, подчас наивных и странных, но иначе сложно научить их думать. Дети должны познавать математику активно, споря между собой, совершая ошибки, за которые никто не поставит двойку, приобретая опыт маленьких открытий. Важна и увлеченность учителя, он должен любить решать задачи, коллекционировать идеи решения и типичные ошибки, создавать новые подборки задач, выстраивать по ходу занятия цепочки задач, уметь давать минимальные подсказки, если задача не решается за разумное время.

Но нередко у молодых преподавателей занятия складываются не совсем удачно. Причин бывает много -– неумение держать дисциплину, плохая обратная связь, слишком быстрый темп, трудные или неинтересные для детей задачи, сложные языковые конструкции.
Могут помешать и однообразные формы работы.
Так, в ВМШ (Вечерняя математическая школа) было важно вовлечь всех детей в решение задачи. Для этого подходят разные игровые формы: голосование за различные ответы, личные и командные соревнования, учет личных достижений, обсуждение парадоксов, исторические байки и т.д.

Диалог хорош тем, что он компенсирует недостаток школьной методики обучения, когда используется только индивидуальная форма работы, а помогать друг другу не положено (подсказка, списывание) и коллективное решение задач не применяется, хотя научиться работать в коллективе тоже важно. Но дети интуитивно чувствуют полезность
коллективной работы и охотно на нее отзываются. При этом учитель ``дирижирует'', подкидывает ``дровишки'', согласует
темп, упорядочивает работу, подводит итоги. Особенно приятно, когда рождается что-то новое, чего учитель не ожидал (идея, задача, шутка и т.д.). Достаточно удачной формой проведения занятия является лекция с диалогом о ходу (взаимные вопросы, голосование и пр.).

\subsection{Открытые вопросы}

Коллективное обсуждение особенно важно при исследовании ``недетерминированных’’ задач (задач ``с открытым ответом''). К сожалению,
почти
все задачи четко формулируются -- что дано и что надо получить (доказать). Это в определенной мере наводит на путь решения и не способствует развитию
творческого поиска. (Когда мы работаем с глубокой задачей на занятии, мы иногда можем давать её как недетерминированную, слегка поменяв вопрос, но тогда мы выходим за рамки листковой системы. Да и таких задач становится все меньше \cite{dver}.)
\footnote{
Как сообщает В.~Ю.~Губарев, в Новосибирске под руководством С~.В.~Августиновича проводится обсуждение интересных задач, придумывание новых задач, обсуждение возможных обобщений, других формулировок,
с чем это может быть связано в математике и т.д. Так появляются и темы для исследования школьниками. Например, при обсуждении пропорций, которые можно получить с бидонами
ёмкостью 1 л и 3 л,
возникла следующая задача (Новосибирск, школьный тур Всероссийской математической олимпиады, 2013 г., 7--8 кл.). {\it Есть две цистерны с неограниченными запасами кофе и сливок (в одной кофе, в другой сливки) и цистерна, в которую можно сливать неограниченное количество жидкости. Есть также два бидона вместимостью $1$ и $3$ литра. Как с их помощью получить
$1$ литр напитка, $5/12$ которого составляет кофе, а $7/12$~-- сливки?}}

\subsection{ Вбрасывание задач}
Это способ организации диалога. Так же как и в системе листков,
в основе стоит задачный подход. Однако здесь мы опираемся на публичность.

Учитель вбрасывает задачу и дает несколько минут на размышление. Если задача несложная, то учитель ходит по рядам и проверяет решения. Первый решивший получает право быть ассистентом учителя, т.е. проверять решения этой задачи у других. Если задача трудная, то существует форма коллективного решения -- ученики набрасывают идеи решения, а учитель отбирает из них перспективные. На этом пути могут понадобиться минимальные подсказки, причем важно спросить детей: ``Ну, что, подсказать немножко?'' Обычно они отказываются.

\subsection{Игровые формы работы}

Игровые формы работы полезны для младших возрастов, особенно в том случае, когда группа устала. Игровая форма должна быть завершающей -- после этого школьников надо распустить. Математические игры также снимают многие недостатки листковой системы. Особенно важно слушать друг друга на математических боях (см.\cite{Matboi}) .

Бывают задачи с подвохом, в которых легко ошибиться, если пытаться угадать ответ или искать его подбором. Учитель не проверяет решения, а выписывает разные ответы на доске. Через 2--3 минуты проводится голосование. Бывает поучительно, когда правильный ответ не набирает большинства голосов. Важно, чтобы работал весь класс, каждый в меру своих сил. Большую роль при этом играет эмоциональная составляющая, поэтому учитель может пошутить, рассказать историю из жизни великого математика, выделить интересные идеи решения, кого-то публично похвалить. Выступая, дети учатся говорить, корректно спорить, выслушивать оппонента. Хорошо, когда атмосфера непринужденная, когда дети учатся не для отметки, а для себя, когда идет поиск истины.

Приведем сюжет, где полезно голосование.

\subsubsection{Можно или нельзя?}

Целью этого занятия является пропедевтика  {\it математического доказательства}, мы демонстрируем на примерах, что, казалось бы, ``очевидные'' вещи
могут оказаться
неверными, формируем потребность в строгих математических рассуждениях. Обсуждение начинается со следующего  примера:

\medskip
{\bf Серия 1.\ 1.} {\it Из трех палочек не всегда можно составить треугольник. А из $100$ палочек
всегда ли
можно найти три, из которых треугольник составляется? }
\medskip

Этот вопрос голосуется, и большинство детей голосуют -- ДА, всегда составляется. Затем кто-то находит контрпример (степени двойки и т.п.).

Хорошо, продолжаем, {\bf 2.} {\it пусть у нас есть $12$ палочек, из которых уже сложены $4$ треугольника. Верно ли, что из них всегда можно сложить $3$ четырехугольника?}

Наученные горьким опытом, дети голосуют ``нет, не всегда''. {\it Мы установили истину, или контрпример все же требуется?}~-- говорит преподаватель. Строятся контрпримеры, но все не подходят. Всякий раз треугольник с самым маленьким периметром
можно ``разобрать на запчасти'' -- вставить его стороны в остальные треугольники. В итоге получается доказательство того, что всегда можно.

Теперь переворачиваем задачу. {\bf 3.}\ {\it Пусть у нас уже есть $3$ четырехугольника. Можно ли из них составить $4$ треугольника?}

Опять большинство голосует за то, что всегда можно, затем возникает контрпример, после чего обсуждаем вопрос о том, а {\it всегда ли можно составить хотя бы один треугольник?} И это, оказывается, не всегда!

Далее по аналогичному сценарию разбирается тройка задач:

\medskip
{\bf Серия 2.}\
{\bf 1.\ } {\it  В коробке лежат карандаши. Есть два карандаша разного цвета и два карандаша разного размера. Верно ли, что найдутся два карандаша, различающиеся и по цвету, и по размеру?}

\medskip
{\bf 2.\ } {\it Девочки в классе различаются ростом, весом и размером обуви. Верно ли, что найдутся две девочки, отличающиеся всеми тремя параметрами одновременно?}

\medskip
{\bf 3.\ } {\it В коробке лежат карандаши. Есть карандаши трех разных цветов и трех разных размеров. Верно ли, что найдутся три карандаша, попарно различающиеся и по цвету, и по размеру?}
\medskip

Обсуждаются и некоторые арифметические задачи. Например:
\medskip

{\bf 4.\ } {\it Группа граждан страны А эмигрировала в страну Б.
Может ли средний IQ жителей обеих стран возрасти?}
\medskip

В какой-то момент школьники сердятся: ``А может, мы сперва подумаем, а затем голосовать будем?''.

\section{Примеры тем и целей занятий, для которых листковая система недостаточна}

Мы приведем несколько примеров занятий (для учащихся самого разного уровня -- от начинающих до продвинутых, вплоть до студентов), при проведении которых необходим диалог с классом. Он нужен, в частности, при воспитании {\it теоретического мышления}. Конечно, приводимые ниже комментарии можно внести в листок. Но дело в том, что подход к математике, как науке о решении занимательных задач и головоломок, является привычным и его разъяснять не надо. В его рамках проще организовать самостоятельную деятельность.
Иное дело, когда ученик встречается (сейчас, увы, реже, чем раньше) с принципиально
другим
отношением к математике, например, когда объясняется, почему конструкции и понятия именно такие, а не другие, и как они строятся. Здесь необходим диалог.

В 7-м классе, например, детям трудно придумать определение треугольника, они могут дать множество неверных определений. В этой ситуации хороша игра: «Кто опровергнет данное определение?», например, ученик говорит: ``это три точки и соединяющие их отрезки, которые образуют три угла''. Всему классу кажется, что определение верное, а учитель рисует на доске три отрезка, выходящие из одной точки и образующие три угла… И не надо жалеть времени на игру в определения -- ученикам полезно ``выстрадать'' определение треугольника как замкнутой трехзвенной ломаной, они поймут цену строгим определениям, приобретут вкус к четкому мышлению.

\subsection{Проблемная ситуация и ее модель}     \label{SbScProblSit}

В школе и вузе практически не учат постановкам задач, а это насущная ситуация и в фундаментальной науке, и в прикладных исследованиях. Приведем пример задачи, в которой важно осознать
чёткую постановку.

\medskip
{\bf Задача 1.} {\it Солнце в зените, над плоской площадкой висит вертолет. Вопрос: тень от вертолета больше него, меньше него или равна вертолету?}
\medskip

Обычно в классе появляются сторонники всех трех гипотез.

{\bf 1.} {\it Больше вертолета}, потому что солнце далеко и его можно считать точкой, тогда получаем расходящийся конус лучей (это видно, когда лучи пробиваются между тучами).

{\bf 2.} {\it Равна вертолету}, потому что от Солнца идет параллельный поток лучей (золотой дождь), и под ним все размеры сохраняются.

{\bf 3.} {\it Меньше вертолета}, поскольку Солнце больше вертолета, и от него к вертолету идет сходящийся конус лучей.

При голосовании обычно максимальное число голосов набирает первая гипотеза. Затем следует обсуждение. В 6--7 классе оно бывает столь бурным, что учитель получает уникальный опыт наблюдения за своими учениками.

Самое интересное, что если в классе нет грамотного астронома, то еще не было случая, чтобы ученикам удалось договориться. В какой-то момент учитель вмешивается в ситуацию и задает неожиданный вопрос: {\it заметил ли кто-нибудь, что в задаче не дано определение тени? Так что же такое тень? Какие точки плоскости мы будем называть тенью?} Следует замешательство, но обычно кто-то говорит, что {\it тень --- это откуда не видно солнце}. Да, -- говорит учитель, -- это одно из определений {\it тени}, но можно дать и другое. Вот мы говорим ``тень от дерева’’, но ведь между листьями солнечные лучи пробиваются --- значит, это не полная тень. И появляется второе определение: {\it Тень --- это точки, в которых Солнце
хотя бы частично
загорожено}. В астрономии говорят, что это область {\it полутени}.

Теперь мы можем понять, что часть людей понимает под {\it тенью} полную тень, а часть –- полутень. Рисуется картинка, на которой ясно, что полная тень меньше вертолета, а полутень -- больше. Получается, что ответ зависит от определения тени. Это и есть пробел в постановке задачи, который никто не заметил.

\medskip
{\bf Задача 2.} {\it В единичный квадрат бросили $101$ точку,
причём никакие три не лежат на одной прямой. Докажите, что найдется треугольник с вершинами в этих точках, площадь которого не превосходит $1/100$.}
\medskip

Задача решалась на математическом бое 10-х классов, и один из авторов сидел в жюри, происходила проверка корректности вызова. Неожиданно для жюри докладчик рассказал удивительно простое и короткое решение. Вот оно:

Выберем произвольную точку и соединим ее со всеми остальными, -- получится $100$ отрезков. Выберем направление по часовой стрелке и последовательно соединим концы отрезков -- получится $100$ непересекающихся треугольников, суммарная площадь которых не превосходит $1$. По принципу Дирихле, найдется треугольник площади не больше $1/100$.

С тех пор эта красивая ошибка стала предметом разбора на занятиях, а ее анализ – поводом для разговора о том, что такое строгое доказательство. \footnote{Если все точки образуют вершины выпуклого многоугольника, то получится $99$ непересекающихся треугольников, и решение не проходит.}

Другие аналогичные примеры содержатся в статье \cite{GeinKovaldjiSapir}. В книге \cite{Akvarium} обсуждается забавный сюжет: {\it Сколько раз надо перегнуть газету, чтобы она достала до Венеры?}. Вычисления дают ответ -– 50 раз. Даем школьникам газету и просим ее перегнуть. По сути дела, во многих суждениях присутствует {\it неявное} или {\it скрытое моделирование}. Оно возникает и тогда, когда выписывается та или иная математическая формула. О преподавании математики нематематикам см. \cite{BelovSnider}.

\subsection{Построение теории по аналогии. Параллельные миры}

В математике довольно часто встречаются аналогичные ситуации и объекты. В этом случае теория строится по образцу теории для аналогичного объекта, переносятся идеи и конструкции, обычно в обе стороны. Мы приводим две темы занятий, содержанием которых является {\it не фактический материал сам по себе} (его несложно дать в виде листков), а именно {\it развитие теории по аналогии}. Параллелизм между мирами необходимо обсуждать, это достаточно непривычно для школьников. Здесь необходим диалог преподавателя с учащимися.

\subsubsection{Геометрия остатков}

Занятие начинается с вопроса: {\it сколько ферзей можно поставить на шахматной доске $5\times 5$ так, чтобы они не били друг друга?} Ответ -- 5, и расставляются они ходом коня с переносом, т.е. доска рассматривается как торическая. Аналогично действуем для доски $7\times 7$. Теперь поставим аналогичный вопрос для магарадж ({\it магараджа}~-- это фигура, которая ходит и как ферзь, и как конь). Попробуем расставить 7 магарадж на доске $7\times 7$, ставя их ходом большого коня $(3,1)$, и получаем неудачу. А вот на доске $11\times 11$ все получается.

Обобщим задачу. Назовем {\it $k$-монстром} фигуру, которая ходит
вдоль прямых $ax+by=0, |a|, |b|\le k$. Ферзь -- это $1$-монстр, а магараджа -- $2$-монстр. Оказывается, при простом $p\ge k^2+1$ на шахматной доске $p\times p$ можно расставить $k$-монстров в количестве $p$ так, чтобы они не били друг друга.

Попробуем это доказать. Вернемся к расположению ферзей на доске $5\times 5$. Запишем их координаты: $\{(0,0),(1,2),(2,4), (3,1), (4,3)\}$ и постараемся дать {\it простое описание} этого множества. Начнем нумерацию с нуля. Получаем $y\equiv 2x \ \ (\!\!\mod 5)$. Соответственно, $k$-монстры ставятся вдоль прямой $y=(k+1)x \ \ (\!\!\mod p)$. Если при этом два монстра будут бить друг друга, то линия боя монстров пересечется с линией, на которую их поставили, хотя бы дважды. Что невозможно: две прямые, как на дискретной, так и на обычной плоскости, пересекаются не более чем в одной точке. Задача про монстров решена.

Вернемся к задаче 2 из раздела \ref{SbScProblSit}.

\medskip
{\it В единичном квадрате расположена $101$ точка, никакие три не лежат на одной прямой. Докажите, что найдется треугольник площади меньше $1/100$ с вершинами в отмеченных точках.}
\medskip

Квадрат режется на $50$ равных полосок, в одну из которых попадет
хотя бы $3$
точки. Любой треугольник внутри прямоугольника занимает не больше половины площади. Последнее утверждение поучительно доказывается с помощью линейного варьирования: ,удем менять треугольник, чтобы его площадь возрастала. Если есть вершина, которая не лежит на границе прямоугольника, то
передвинем
ее так, чтобы она попала на границу. Если есть вершина треугольника, которая не лежит в вершине прямоугольника, то передвинем ее в одну из соседних вершин так, чтобы площадь увеличилась. Теперь площадь треугольника равна половине площади прямоугольника. Этот метод требует обсуждения, о нем сказано в \cite{KnelKovalKatastrofy}.

Затем переходим к {\it нижним оценкам}. Найдём такое расположение точек,
 что площадь любого треугольника с вершинами в них не меньше $\frac{1}{2\cdot 101^2}$. Увеличив единичный квадрат в $101$ раз, приходим к задаче:

\medskip
{\it В квадрате $101\times 101$ расположить $101$ целочисленную точку так, чтобы никакие $3$ не попали на одну прямую.}
\medskip

Строим естественную конструкцию. Первая точка имеет координаты $(0,0)$, вторая $(1,0)$, далее $(2,1)$, $(3,3)$ и т.д. На $k$-м шаге новая точка смещается на 1 вправо и поднимается на $k$ вверх, по сравнению с предыдущей. Как и в случае расположения монстров, когда достигнута верхняя сторона квадрата, переходим на нижнюю по циклу (считаем квадрат тором). Но почему никакие три из построенных точек не попадут на одну прямую? Как и в задаче про монстры, запишем их координаты $(x,y)$. Оказывается, $y\equiv x(x-1)/2 \ \ (\!\!\!\!\mod 101)$. И то, что никакие три из отмеченных точек не попадают на одну прямую, отражает тот факт, что парабола с прямой не могут пересечься в трех точках.

Далее обсуждаются магические квадраты. ({\it Магический квадрат}~-- это расположение чисел от 0 до $n^2-1$ в квадрате $n\times n$ такое, что суммы по вертикалям, горизонталям и главным диагоналям равны.) Сперва в точку с координатами $(k,l)$ ставим число $nk+l$ (нумерация координат идет от нуля). Если теперь выбрать {\it ладейное множество} полей (по одной клетке на вертикали и по одной по горизонтали), то суммы
в любом таком множестве будет одна и та же.
Любая прямая, кроме горизонтали или вертикали, есть ладейное множество. Итак, вдоль любой прямой, кроме тех,
которые нам нужны (вертикалей и горизонталей), суммы постоянны.

Теперь мы выберем две непараллельные прямые (не вертикали и не горизонтали) $L$ и $M$. Рассмотрим их как {\it оси координат}. Перерисуем таблицу. В новой таблице в точке с горизонтальной координатой $x$ и вертикальной координатой $y$ будет стоять такое же число, как в старой таблице, с теми же координатами, но в системе $(L,M)$. В новой таблице горизонтали отвечают расположениям чисел вдоль прямой, параллельной $L$, а вертикали -- вдоль прямой, параллельной $M$. Поэтому суммы по горизонталям и вертикалям равны. (Еще остается позаботиться о том, чтобы при получившемся преобразовании старая вертикаль или горизонталь не оказалась диагональю).
Получаем искомый магический квадрат.

В процессе занятия напоминаем, что такое косоугольная система координат на плоскости, и переносим конструкцию на дискретную плоскость -- в данном случае на доску $5\times 5$. Рисуем ``косые'' оси и числа, в них стоящие, заполняем новую таблицу ``руками''.

Возникающая отсюда задача позволяет обсудить понятие конечной проективной плоскости:

\medskip
{\it Городок DIV-GRAD обладает отличной автобусной сетью! С каждой
остановки можно проехать на любую другую без пересадок.
Каждый маршрут имеет пять остановок, а каждые два маршрута
имеют единственную общую остановку. Сколько же автобусных маршрутов
в этом дивном граде? Нарисуйте карту.}
\medskip

Таким образом, это занятие демонстрирует аналогию между обычной геометрией и дискретным миром и их взаимосвязь.
\footnote{Для продвинутых учащихся полезно обсудить следующую задачу (\cite{kvadrat}):
{\it Квадрат разбит на треугольники равной площади. Доказать, что их число четно.}

Идея доказательства:
Если количество треугольников нечетно, то  удвоенная площадь каждого сравнима с нулем по модулю $2$, так что векторное произведение векторов, образующих треугольник, равно нулю по модулю два. Тем самым любые три вершины треугольничков разбиения коллинеарны по модулю $2$ (если координаты вершин треугольников разбиения есть рациональные числа, в знаменатель которых не входит двойка),  а вершины единичного квадрата -- нет. Общий случай отличается добавлением алгебраической техники.  }

При ведении занятий по этой теме, чтобы  подчеркнуть параллелизм между арифметическим и вещественным миром, рисуется доска, разграфленная на две колонки, например:

\begin{tabular}{|l|c|}
 \hline
 Вещественные числа & Остатки по модулю $p$ \\
 Вещественная плоскость & Клетчатая доска \\
 Прямая на плоскости & Линия боя монстра, линия расстановки монстров \\
 Две прямые пересекаются не более чем в одной точке & Монстры не бьют друг друга \\ Косоугольная система координат & Перерисовка квадрата \\
  Проективная плоскость & город DIV-GRAD \\
 \hline
\end{tabular}

\subsubsection{Двойные и комплексные числа, ортогональность и псевдоортогональность}  \label{DobleNumbers}

Рассмотрим числа вида $a+bj$, где $j^2=+1$. Для них строится теория, параллельная теории комплексных чисел. Листок осмыслен только как дополнение к рассказу. Цель занятий -- показать {\it параллелизм} теорий и построение теорий по аналогии. Прежде всего, надо упредить недоразумение -- если $j^2=1$, то школьники могут считать, что $j=\pm 1$, а это –- не так, $j$ есть формальная буква. В свете этого надо обсудить формальное определение комплексного числа как пары вещественных чисел и аналогичное определение {\it двойного числа} (вместо термина ``двойное число'' удобнее употреблять термин ``гиперболическое'') $a+bj$ опять-таки, как пары вещественных чисел, но с несколько другим умножением (в случае комплексных чисел $(a,b)\cdot (c,d)=(ac-bd,ad+bc)$, в гиперболическом случае $(a,b)\cdot (c,d)=(ac+bd,ad+bc)$). Далее, легко определить квадрат модуля как $a^2-b^2$ и проверить свойство произведения. Гораздо важнее понять, почему определение таково.

 Рассмотрим аналог комплексных чисел с условием
 $i^2=-c,\ c>0.$
 Тогда, положив $i'=i/\sqrt{c}$, имеем $(i')^2=-1$, и мы возвращаемся к теории обычных комплексных
чисел. При этом $|z|^2=a^2+cb^2$. Теперь положим $c=-1$ и получим, что $|a+bj|^2=a^2- b^2 $. Аналогом единичной окружности становится гипербола $x^2-y^2=1$. Далее обсуждается концепция {\it аргумента} как удвоенной площади соответствующего криволинейного сектора. По аналогии с ситуацией, когда $i^2=-c$, показывается, что при умножении аргументы складываются.

В заключении можно поговорить о {\it дуальных числах}, когда $j^2=0$.

Главное здесь -- рассказ о {\it построении теории} по аналогии, в данном случае речь идет о приеме введения параметра и подстановки значения из ``запретной области'',
в данном случае из области $c<0$.

\subsubsection{О геометрии Лобачевского}
Тема, изложенная в разделе \ref{DobleNumbers} служит пропедевтикой к изучению неевклидовой геометрии (см. \cite{YaglomGalilei}).
Одна из лучших ознакомительных книг -- \cite{Shirkov}. Отметим только псевдонаучность следующего, увы, распространенного подхода: описывается модель (Кэли-Клейна или Пуанкаре) и проверяются аксиомы. Такой листочек легко составить, но смысла в таком занятии (тем более, не поддержанном обсуждением) будет немного. При такой подаче модель не возникает естественным образом. Модель Кэли-Клейна естественно возникает, например, при придании смысла понятию {\it сфера мнимого радиуса} и построении ее центральной проекции.
\footnote{В истории неевклидовой геометрии есть загадка, на которую авторы не имеют удовлетворительного ответа. Оказывается, что и Лобачевский, и Бойяи, и Гаусс имели модель Кэли-Клейна, но при этом непротиворечивость неевклидовой геометрии так ими и не была доказана! Они рассматривали плоскость, касательную к орисфере, и ее проекцию на орисферу вдоль пучка прямых, соответствующего орисфере. Внутренняя геометрия орисферы -- обычная евклидова геометрия. Неевклидова плоскость проецировалась на внутренность круга, а прямые -- на хорды. Все формулы
были известны (см. \cite{Shirkov}). Оставалось только сказать слова ``модель'' и ``назовем {\it плоскостью} внутренность круга, {\it прямой}~-- хорду, {\it расстоянием}~-- логарифм двойного отношения четверки точек, возникающих при сечении'' (формула, известная Лобачевскому, Бойяи и Гауссу). Возможно, они попали в ловушку кантовской философии: геометрические понятия даны Богом, отсюда их истинность. Нам представляется разговор об этом со школьниками
чрезвычайно важным. При открытии специальной теории относительности также пришлось преодолевать кантовскую философию и ситуация была во многом аналогична открытию неевклидовой геометрии. Об этом пишет Вертгеймер \cite{Wertgeimer} (создатель гештальт-психологии, по совместительству учитель математики, он интервьюировал Эйнштейна). В его книге обсуждаются интересные задачи, которые можно использовать для занятий.}

 Один из авторов сталкивался с ситуацией, которая показательна
 как тип поведения тренера. Один участник сборов на международную олимпиаду хорошо умел решать задачи, используя {\it поляритет}
\footnote {См. https://ru.wikipedia.org/wiki/Полюс\_и\_поляра},
 который воспринимался как олимпиадный трюк. О связи поляритета с псевдоскалярным произведением и соответствием ``экватор--полюс’’ на сфере школьники не подозревали. На вопрос тренеру о том, почему школьникам это не рассказывали, был получен ответ ``нет времени'', т.е. смысл понятия и его научное содержание тренеру неважно.

\subsection{Ознакомление с новыми
концептуально значимыми понятиями}

Мы приведем несколько примеров работы с новыми понятиями, где необходим диалог.

\subsubsection{Метод математической индукции}

Формирование понятия математической индукции -- дело далеко не простое. Соответствующие методические проблемы ставятся и решаются в замечательной работе И.~С.~Рубанова \cite{Rubanov}. В  ней подробно разбираются, в частности, диалоги, возникающие при преподавании этой темы. Методические проблемы формирования новых понятий обсуждаются также в работе \cite{BelovYavich}, которая основана на статье И.~С.~Рубанова. Работа \cite{Rubanov} интересна и важна отнюдь не только в плане преподавания специфической темы, но и в более общем контексте.

\subsubsection{Теория информации}

 Цель занятия -- не только
обучение некоторым приемам,
предназначенным для решения олимпиадных задач, но и формирование концепции количества информации. Листок может дополнять, но не заменять диалог.
Эту тему можно излагать для учащихся разного уровня
-- от начинающих до продвинутых. Для начинающих следует ограничиться задачами о взвешивании и угадывании (см.\cite{MihalinNikonov,BelovFrenkinMihalinNikonov}), объяснением, что такое {\it трит} и как он соотносится с {\it битом} (понятие {\it логарифма} здесь можно объяснить по
ходу дела). Для более продвинутых можно затронуть вероятностные аспекты. И, наконец, для наиболее продвинутых можно затронуть понятие {\it энтропии} (см. \cite{Yaglom}). В последние годы эта тема (кроме задач на взвешивание) в кружках практически не представлена, но в то же время, в связи с развитием информатики, ее значение возросло. Поэтому мы ее разбираем подробно.

Занятие начинается со вступления. {\it Все мы сталкивались с понятием ``информация'' и ее измерением. Мы знаем, что на флешке помещается меньше информации, чем на жестком диске. Но как выразить математически это бытовое понятие? Например, утверждается, что 15 июля в Сахаре не будет дождя. Такая фраза вызывает улыбку. Тем не менее, какая-то информация заключена в этом прогнозе. Какая именно, и как это измерить?}

Мы обсуждаем задачи об угадывании задуманного числа, меньшего 1000, за минимальное число вопросов, на каждый из которых возможен ответ ``да'' или ``нет''. Далее обсуждается вопрос об обнаружении более тяжелой фальшивой монеты на чашечных весах без гирь за минимальное число взвешиваний. В первом случае среди $2^n$ чисел за $n$ вопросов можно угадать задуманное (а из большего набора нельзя), а во втором -- из $3^n$ монет можно определить фальшивую. При этом педалируется идея {\it пространства вариантов} и его сжатия после каждого действия.

Далее обсуждается случай, когда неизвестно, легче ли фальшивая монета или тяжелее и как ее определить из 12 монет за 3 взвешивания, или из 39 монет -- за 4 взвешивания. Здесь идея пространства вариантов помогает сократить перебор -- взвешивания надо подбирать так, чтобы число возможностей при
максимальном исходе было минимальным. Подробности см. \cite{MihalinNikonov,BelovFrenkinMihalinNikonov}.

Теперь мы говорим о {\it битах} и {\it байтах} и о {\it тритах} и объясняем, почему трит -- это $\log_23$ бит.
(Винчестер с $k$ битами может находиться в $2^k$  состояниях, а винчестер с $l$ тритами может находиться в $3^l$ состояниях. Информацию с одного можно перегнать на другой и обратно, если $2^k\sim 3^l$. Отсюда получается, что $k\sim \log_23\cdot l$, т.е. $\log_23$ бита в пересчете на один трит при $k\gg 1$.)
Аналогично вводится и обсуждается {\it $n$-ит}, содержащий $\log_2n$ битов.
Здесь можно попутно обсудить с учащимися понятие {\it логарифма}.

Для более продвинутых учащихся можно углубиться далее. Возвращаемся к угадыванию чисел. Если задумано число, меньшее 1000, и мы задали вопрос: ``задуманное число больше 500?'', то мы получаем один бит информации из ответа. Хорошо, а если зададим вопрос ``задуманное число больше 100?'', то какую информацию мы получим?

С вероятностью $0.1$ мы сократим количество возможностей до $100$ и получим информацию, равную $\log_21000-\log_2100=\log_210$,
Но с вероятностью $0.9$ мы сократим количество возможностей только до $900$ и получим информацию, равную $\log_21000-\log_2900=\log_2(10/9)$. Следовательно, матожидание полученной информации будет равно $0.1\cdot \log_210+0.9\cdot \log_2(10/9)$. Аналогично, если есть опыт с возможными исходами, вероятности которых равны $p$ и $q$ ($p+q=1$), то матожидание количества информации будет равно $-p\log_2p- q\log_2q $,
а если возможных исходов несколько и каждый из них имеет вероятность $p_i$ ($\sum p_i=1$), то ожидание количества информации будет равно $-\sum p_i\log_2p_i$.

Здесь мы уже готовы формально вычислить информацию, которую несет утверждение о том, что 15 июля 2015 года в Сахаре будет ясная погода (вероятность дождя $10^{-4}$), но тем не менее надо понять его {\it смысл}. Рассмотрим последовательность показаний некоторого датчика.
Вероятность регистрации сигнала мала, скажем, $10^{-4}$, но последовательность показаний достаточно длинная, длины $n\gg 1$. Итак, мы имеем длинную последовательность нулей и единиц, с количеством единиц примерно равным $10^{-4}n$. Двоичный логарифм от числа таких последовательностей и есть информация, которую выдал нам датчик. Разделив эту величину на $n$ и перейдя к пределу, мы поймем, что означает информация, которую дает одно показание датчика.

Материал такого рода может варьироваться в зависимости от уровня (силы) класса. С достаточно сильными учащимися можно обсудить энтропию и теорию Шеннона.

Все вышеперечисленные вопросы нуждаются в публичном обсуждении перед классом, ибо целью является формирование новых понятий и
методологии. Листки здесь могут быть полезны только как дополнение.

\subsection{От решения открытой проблемы - к методу}

\subsubsection{Квадратичный закон взаимности и построение правильных $n$-угольников}

Данная тема предназначена для достаточно продвинутых учащихся.

К.~Ф.~Гаусс одновременно решил две задачи (публикации различались по времени на неделю). Он выяснил, при каких простых $p$ в поле вычетов ${\mathbb Z}_p$
существует $\sqrt{q}$ и одновременно построил правильный семнадцатиугольник. Как он действовал?

Легко видеть, что если $\sqrt{q_1}$ и $\sqrt{q_2}$ принадлежат ${\mathbb Z}_p$, то $\sqrt{q_1q_2}\in{\mathbb Z}_p$, если же $\sqrt{q_1}\in{\mathbb Z}_p$ и $\sqrt{q_2}\notin{\mathbb Z}_p$, то $\sqrt{q_1q_2}\notin{\mathbb Z}_p$. Несколько труднее показать, что если $\sqrt{q_1}\notin{\mathbb Z}_p$ и $\sqrt{q_2}\notin{\mathbb Z}_p$, то $\sqrt{q_1q_2}\in{\mathbb Z}_p$.
Тем самым задача сводится к выяснению вопроса о принадлежности ${\mathbb Z}_p$ величин 
$i=\sqrt{-1}=\sqrt[4]{1}$ и $\sqrt{q}$ при простых $q$.

Используя малую теорему Ферма (в обратном направлении) и тот факт, что многочлен  $n$-й степени имеет
не больше $n$ корней, получаем, что $i\in \{0,1,\dots,p-1\}$ тогда и только тогда, когда $i^p-i=0$. Это -- {\it первая ключевая идея}, нуждающаяся в акцентировании внимания и обсуждении.

Естественно получить более общий факт: если $\xi=\sqrt[n]{1}$~-- неприводимый корень $n$-й степени из единицы, то $\xi^p-\xi=0 \iff p\equiv 1 \ \ (\!\!\mod n)$. В частности, $\frac{-1+\sqrt{-3}}{2}=\sqrt[3]{1}\in{\mathbb Z}_p \iff p=3k+1$. Следовательно, $\sqrt{-3}\in{\mathbb Z}_p \iff p=3k+1$.

Окрыленные этим успехом, получаем, что
$\sqrt{2}=\omega+\omega^{-1}$, где $\omega=\sqrt[8]{1}$. Поэтому
 $\sqrt{2}^p=(\omega+\omega^{-1})^p=\omega^p+\omega^{-p}=\omega+\omega^{-1}=\sqrt{2}$ так что
 $\sqrt{2}^p=\sqrt{2}$   при $p\equiv\pm 1 \ \ (\!\!\mod 8)$ и только в этом случае.
  Теперь ясно, что $\sqrt{5}$ надо выразить через корни пятой степени из 1, т.е. через корни уравнения: $x^4+x^3+x^2+x+1=0$. Разделив на $x^2$ и сделав замену $t=x+1/x$, находим $t=(-1\pm\sqrt{5})/2$. Итак, если  $\xi=\sqrt[5]{1}$~-- неприводим,
 то $t=\xi+\xi^{-1}$.
Рассуждая аналогично,
получаем, что $\sqrt{5}\in{\mathbb Z}_p$, если $p\equiv \pm 1 \ \ (\!\!\mod 5)$.

По аналогии с предыдущими случаями, стремимся выразить $(-1\pm\sqrt{\pm 7})/2$ через корни уравнения $x^6+x^5+x^4+x^3+x^2+x+1=0$. Наивное копирование предыдущего случая не проходит: если разделить на $x^3$ и сделать замену $t=x+x^{-1}$, то возникает уравнение третьей степени. И здесь нужна

{\it Вторая ключевая идея}. Что такое неприводимый корень седьмой степени из единицы? Это буква, удовлетворяющая равенству $x^6+x^5+x^4+x^3+x^2+x+1=0$. С этой точки зрения $x$ и $x^k$ ничем не отличаются (если только $k\not\vdots 7$). Когда корни разбили на пары ($x+x^{-1}$), их
получилсь
три. Замены типа $x^k\to x$ эти три пары  переставляли, не нарушая группировки.
При этом суммы $x^k+x^{-k}$ оказались корнями многочлена третьей степени. Значит, надо {\it сгруппировать корни не в три, а в две группы, так чтобы замена вида $x^k\to x$ сохраняла группировку}. Тогда следует ожидать (по аналогии с предыдущими случаями), что сумма в одной группе равна $(-1\pm\sqrt{-7})/2$. Если $\xi$~-- корень уравнения  $x^6+x^5+x^4+x^3+x^2+x+1=0$, то с $\xi$ группируются
$\xi^2, \xi^4$.
И действительно, легко проверить, что $\xi+\xi^4+\xi^2=(-1\pm\sqrt{-7})/2$.
И если $p\equiv 1,2,4 \ \ (\!\!\mod 7)$, то $\sqrt{-7}\in{\mathbb Z}_p$, если же $p\equiv 3,5,6 \ \ (\!\!\mod 7)$, то $\sqrt{-7}\notin{\mathbb Z}_p$, поскольку тогда $(\xi+\xi^4+\xi^2)^p=\xi^3+\xi^5+\xi^6=1-\xi+\xi^4+\xi^2$. (Если $\xi+\xi^4+\xi^2=1-\xi+\xi^4+\xi^2$, то $(-1\pm\sqrt{-7})/2=\xi+\xi^4+\xi^2=-1/2$,
так что $\sqrt{-7}\equiv 0 \ \ (\!\!\mod(p))$.)

Далее разбиваем корни 11-й, 13-й, 17-й, 19-й степени из единицы на две группы. Соответствующие суммы оказываются равными
$(-1\pm\sqrt{-11})/2$, и т.д.

Теперь надо дать описание группировки корней на две группы в общем виде. Пусть $\xi=\sqrt[n]{1}$~-- неприводимый корень $n$-й степени из единицы. Что  группируется вместе с $\xi$?
Выписываются группировки и затем обнаруживается что всякий раз берутся квадраты по модулю $p$!


Сумма же членов в группе равна $\frac{-1\pm\sqrt{\varepsilon(p)p}}{2}$, где $\varepsilon(p)=1$ при $p=4k+1$ и $\varepsilon(p)=-1$ при $p=4k-1$.

Остается проверить равенство $\sum_{k\equiv s^2 \ \ (\!\!\mod(p))}\xi^k=\frac{-1\pm\sqrt{\varepsilon(p)}}{2}$.

Далее из этого выводится {\it квадратичный закон взаимности:} $\sqrt{p}\in{\mathbb Z}_q\iff \sqrt{q}\in{\mathbb Z}_p$, когда хотя бы одно из простых чисел $p,q$ имеет вид $4k+1$; если же они оба имеют вид $4k-1$, то $\sqrt{p}\in{\mathbb Z}_q\iff \sqrt{q}\notin{\mathbb Z}_p$.

{\bf Вторая основная идея} позволяет строить правильный семнадцатиугольник. Надо сперва сгруппировать 16 неприводимых корней 17-й степени из единицы на две группы по 8, как раньше, затем каждую подразбить на две группы по 4, а каждую группу по 4 -- на две группы по 2 корня. При этом любая подстановка $\xi^s\to\xi\ (s=1,\dots,16)$  должна только переставлять группы, не нарушая их структуры.
Полезная курсовая работа для школьника -- построить правильный 17-и угольник, а также $257$ и $65537$ -угольники с (помощью компьютера!).

Рассказ о группировки корней служит пропедевтикой к теории Галуа, идейно разгружая дальнейшее изучение. Здесь соответствие Галуа (инвариантные подполя) задаются группами корней, и их построение выступает как эмпирический прием.
\footnote{Ловушка в преподавании теории Галуа -- в богатстве идей, содержащихся в коротком тексте, но требующих основательной проработки. (Похожая ситуация и с цепными дробями.)}
Данный рассказ позволяет также продемонстрировать учащимся ``гуманитарный''  способ познания: надо вжиться в героя (в данном случае К.~Ф.~Гаусса) и через это увидеть, как он действовал.

\subsubsection{Задача об $n-2$ треугольниках и линейное варьирование}
Более 100 лет (с 1870 по 1979) стояла открытая проблема:

\medskip
{\it $n$ прямых общего положения (т.е. никакие две не параллельны, никакие три не пересекаются в одной точке) делят плоскость на части. Доказать, что среди частей разбиения найдется не менее $n-2$ треугольников.}
\medskip

Эта проблема и вопрос о том, как додуматься до ее решения, обсуждается в
\cite{KnelKovalKatastrofy}. Из ее решения органично возникает {\it метод катастроф}. Работу \cite{KnelKovalKatastrofy} также можно использовать как материал для занятий по этому методу.

\subsection{Демонстрация силы понятия}

\subsubsection{Метод центра масс}

Маленькие дети любят машинку с не слишком сложным устройством: здесь батарейка, здесь выключатель, здесь моторчик, здесь шестеренка. Небольшую теорию, позволяющую решать задачи методом центра масс, можно уподобить такой машинке. Имеется основное понятие -- {\it центр масс}. Имеются две технические леммы: о том, что подсистему грузов можно заменить на груз суммарной массы, расположенный в центре масс подсистемы, и от этого центр масс системы не изменится, и лемма о положении центра масс системы из двух точек с массами $m$ и $M$ соответственно. Миниатюрная работающая теория служит хорошим {\it образцом} научной теории.

Подобную демонстрацию можно повторить при преподавании тем ``решение задач с помощью момента инерции’’, ``аффинные и проективные преобразования’’.

\subsubsection{Алгебра, теория полей}

Цель данной темы – продемонстрировать алгебраические концепции, возникшие в начале 19 века, изменившие лицо алгебры.

Пусть $K$~-- некоторое {\it поле}, т.е. множество (комплексных) чисел, замкнутое относительно арифметических операций, т.е. операций сложения, вычитания, умножения и деления (кроме деления на нуль). {\it Расширение поля} $K$ набором чисел $\{x_i\}$ есть множество чисел, получающееся из поля $K$ и набора $\{x_i\}$ с помощью арифметических операций. Иными словами, это минимальное по включению поле $K'$, содержащее $K\cup \{x_i\}$.

Эти два понятия выглядят совершенно невинно, однако их достаточно для решения знаменитых задач древности об удвоении куба и трисекции угла. Мы это продемонстрируем, чтобы читатель отнесся с должным уважением и к другим понятиям теории полей, отражающим глубокие идеи, которые необходимо не спеша продумать. Короткие тексты с большой концентрацией идей весьма коварны
(например, при изучении
цепных дробей также сталкиваются с большим числом идей в коротком тексте, отсюда методические проблемы). Их надо изучать {\it очень медленно}.

Доказывая невозможность удвоения куба, полезно обсудить феномен {\it формализации}. Чтобы доказать возможность построения циркулем и линейкой, достаточно понимания процедуры с позиции обычного здравого смысла. Однако для доказательства {\it невозможности} необходимо саму процедуру построения сделать математическим объектом, т.е. {\it формализовать.} (Аналогичные разговоры следует вести при изучении построения одной линейкой.) Прежде всего, если задан единичный отрезок, то несложно проверить, что все длины, которые мы можем построить, выражаются через 1 с помощью арифметических действий и вычисления квадратного корня.

Чтобы разгрузить идейную часть, полезно выдать листок с упражнениями. Можно начать с упражнения, которое доступно читателю, знающему, как делить многочлены с остатком:

\medskip

{\bf Упражнение 1.} {\it Даны многочлен $P$ третьей степени и многочлен $Q$ второй степени, оба с рациональными коэффициентами. Докажите, что если $P$ и $Q$ имеют общий корень, то многочлен $P$ имеет рациональный корень.}

\medskip

Этот факт легко обобщается для произвольного числового поля:

\medskip

{\bf Упражнение 2.} {\it Даны многочлен $P$ третьей степени и многочлен $Q$ второй степени, оба с коэффициентами из поля $K$. Докажите, что если $P$ и $Q$ имеют общий корень, то многочлен $P$ имеет корень в $K$.}
\medskip

{\it Квадратичным расширением $K'$ поля $K$} называется расширение $K$ корнем квадратного трехчлена с коэффициентами из $K$.

\medskip

{\bf Упражнение 3.} {\it Докажите, что каждый элемент из $K'$ либо принадлежит $K$, либо является корнем квадратного уравнения с коэффициентами из $K$.}

\medskip

{\bf Следствие 1. } {\it Пусть поле $K$ есть квадратичное расширение поля $L$. Многочлен $P$ с коэффициентами из $L$ имеет степень 3 и не имеет корня в $L$. Тогда у него нет и корня в $K$.}
\medskip

{\bf Следствие 2.} {\it Пусть поле $K_{i+1}$ есть квадратичное расширение поля $K_i \, i=0,\dots, n$, а многочлен $P$ с коэффициентами из $K_0$ имеет степень 3 и не имеет корня в $K_0$. Тогда у него нет и корня в $K_n$.}
\medskip

При разговоре о доказательстве невозможности удвоения куба с помощью циркуля и линейки (доказательство невозможности трисекции угла аналогично) опять говорим о {\it процедуре формализации}, о том, как формализовать понятие {\it разрешимости уравнения в радикалах}, а также {\it разрешимости уравнения в квадратных радикалах}. Обсуждаем понятие {\it башни (квадратичных) расширений}.

Предположим, что
$\sqrt[3]{2}$
 выражается через квадратные радикалы. Тогда существует цепочка полей $K_0={\mathbb Q}\subset K_1\subset\dots K_n$ таких, что

\begin{itemize}
 \item $\sqrt[3]{2}\in K_n$,
 \item $K_i$ есть квадратичное расширение поля $K_{i-1}$ для любого $1\le i\le n$.
\end{itemize}

Как видим, здесь есть предмет для публичного обсуждения, листок играет вспомогательную роль. Точно так же,
при доказательстве теоремы Руффини-Абеля надо объяснять его ``философию'', а не ограничиваться набором задач.

\section{Синтез листков и диалога}

Сказанное не означает, что форма диалога лучше формы листков, –– эти формы хороши в разных ситуациях: диалог лучше на начальных этапах работы кружка, особенно для 3--5 классов, а по мере развития учеников уместно чаще использовать тематические листки для упражнения в самостоятельной работе или для сдачи зачета.

И в диалоге, и при листках успешно работает система плюсиков, т.е. успехи ученика накапливаются, и он может следить, сколько плюсиков у него набралось. При этом нужна умеренность, чтобы плюсики, как внешний стимул к занятиям, не вытеснили главное –- внутренние стимулы к развитию.

Полезно использовать разработанные листки и в индивидуальной работе, и в диалоге с классом.

Например, можно давать листок на группу, тогда ученики решают его коллективно. Составляются команды по 4--6 человек,
которые
получают одинаковые листки с задачами. Каждый член команды должен
ответить хотя бы одну задачу. На каждую задачу дается три попытки, но отвечать должен один и тот же человек. Либо команда решает, кто какую задачу будет отвечать, либо учитель назначает, кто ответит.

Н.~Х.~Розов писал: ``Есть и другие формы работы со школьниками. Я, например, по ``листкам’’ никогда не занимался –– в кружке Н.~Бахвалова их не было. Но применялось ``монотонное'' изучение теории ``по шагам'' и коллективное обсуждение циклов тематически связанных задач. Тогда не то что ``листков'' не было –– с бумагой было туго, а размножать было просто не на чем. Я всегда предпочитал ``ступенчатый'' метод преподавания (изучение ``подъемом по лестнице'') и даже частично применял его в ходе преподавания на мехмате. По моему мнению, этот метод –– совершенно самостоятельный, весьма эффективно развивающий креативность и достойный специального разбора. А самым выдающимся примером реализации такого метода служит книга И.~М.~Яглома и В.~Г.~Болтянского ``Выпуклые фигуры''. Другой любопытный прием я наблюдал у А.~Н.~Колмогорова. Такая система и сейчас иногда применяется в СУНЦе в ходе исследовательской работы школьников'' \cite{Rosov}.

\medskip
{\bf Благодарности.}
Авторы признательны В.~Ю.~Губареву, Г.~А.~Мерзону, С.~В.~Петухову, Л.~В.~Радзивиловскому, М.~И.~Харитонову, М.~И.~Ягудину за полезные обсуждения. Особая признательность -- Н.~Х.~Розову и И.~С.~Рубанову за плодотворное обсуждение и поддержку. Работа поддержана грантом РФФИ № 14-01-00548.

\end{document}